\documentclass[twoside]{article}
\usepackage{latexsym}

\newcommand{\Section}[1]
  { \section{#1}
    \pagestyle{myheadings} 
    \thispagestyle{myheadings} 
    \markboth{\sc #1}{\sc Blackwell Games}
    \setcounter{equation}{0}}

\renewcommand{\theequation}{\thesection.\arabic{equation}}

\newtheorem{thm}{Theorem}[section]   
\newtheorem{cor}[thm]{Corollary}     
\newtheorem{lem}[thm]{Lemma}         

\newtheorem{definition}[thm]{Definition}   
\newenvironment{defn}{\begin{definition} \em}{ \end{definition}} 

\newtheorem{remark}[thm]{Remark}                 
\newenvironment{rem}{\begin{remark} \em}{ \end{remark}}      
\newtheorem{example}[thm]{Example}
\newenvironment{ex}{\begin{example} \em}{\end{example}}
\newtheorem{notation}{Notation}
\newtheorem{terminology}{Terminology}

\newcommand{\N}{{I\!\! N}}                      
\newcommand{\R}{{I\!\! R}}                      

\newcommand{\cat}{\widehat{\hspace{.1in}}}      

\newcommand{\eop}{\nopagebreak\mbox{ }\newline\hspace*{\fill}$\Box$\newline}   

\newcommand{\val}{\mbox{\rm val}}
\newcommand{\upperval}{\mbox{\rm val}^\uparrow}
\newcommand{\lowerval}{\mbox{\rm val}^\downarrow}
\newcommand{\len}{\mbox{\rm len}}

\newcommand{\prc}{\nopagebreak\noindent\bf Sketch of proof: \rm}
\newcommand{\proof}{\nopagebreak\noindent\bf Proof: \rm}
\newcommand{\expect}[3]{E({#1}\mbox{ \rm vs }{#2}\mbox{ \rm in }{#3})}  
\newcommand{\expectst}[1]{E(\sigma\mbox{ \rm vs }\tau\mbox{ \rm in }\Gamma({#1}) ) }
\newcommand{\expectfin}[2]{E(\sigma\mbox{ \rm vs }\tau\mbox{ \rm in }\Gamma_{#1}({#2}) ) }
\newcommand{\payoffdef}[4]{
{#1}(p) & = & {#3} \mbox{ \rm for }p \in {#2} \\
{#1}(w) & = & {#4} \mbox{ \rm if }w \not\in [{#2}]
}
\newcommand{\heq}[0]{\hspace{.35in}}

\newcommand{\switched}[1]{{#1}_{\mbox{\scriptsize sw}}}
\begin{document}

\title{Blackwell Games}
\author{Marco R. Vervoort 
\thanks{My thanks go to Michiel van Lambalgen and Tonny Hurkens, 
for their guidance, ideas and meticulous proofreading.
This research was partially funded by the
Netherlands Organization for Scientific Research (NWO), under
grant PGS 22-262.
}}
 
\date{\today}

\maketitle
\begin{abstract}
\noindent
Blackwell games are infinite games of imperfect information.
The two players simultaneously make their moves
and are then informed of each other's moves.
Payoff is determined by a Borel measurable function $f$
on the set of possible resulting sequences of moves.
A standard result in Game Theory is that finite games of this type are
determined.
Blackwell proved that infinite games are determined, but only
for the cases where the payoff function is the indicator function
of an open or $G_\delta$ set \cite{Blackwell,Blackwell-2}.
For games of perfect information, determinacy has been proven for games of
arbitrary Borel complexity \cite{G-Borel-1,G-Borel-2,G-Borel-3}.
In this paper I prove the determinacy of Blackwell games over a
$G_{\delta\sigma}$ set, in a manner similar to Davis' proof of determinacy
of games of $G_{\delta\sigma}$ complexity of perfect information \cite{Davis}.

There is also extensive literature about the consequences of assuming AD,
the axiom that {\em all} such games of perfect information are determined
\cite{AD-1,AD-2}. In the final section of this
paper I formulate an analogous axiom for games of imperfect information,
and explore some of the consequences of this axiom.
\end{abstract}
\Section{Introduction}

Imagine two players playing a game of Blind Chess.
The only board they have is in their minds,
and they make their moves merely by announcing them.
Someone who doesn't know the rules
would find a game like this difficult to follow.
If that someone was of a literal bend, he might describe it like this:
\begin{quote}
``There were two players, playing against each other.
The first player said something, and I was told it was her move,
and that she had made the move by saying it.
The other player thought for a while, and then announced his own move.
Then the first player made a move again, then the second player, and so forth.
The moves always sounded similar, something like "pawn from ee-four to ee-five".
So I think they couldn't just say anything, but had to select their moves
 from only a few possible options.
And suddenly they stopped, and shook hands, and I was told that the first
player had won, apparently because of the moves she and her opponent had
played.''
\end{quote}
If no one gave the poor fellow a copy of the rules of Chess,
the way a sequence of moves determines which player wins
would probably seem quite arbitrary.
And our hypothetical observer might be quite impressed
that apparently chess-players are able to memorize this long list
of what the result is of each possible sequence of moves.

Of course, the game of Chess is not really that arbitrary,
and those of us who play chess only need to know a few simple rules
to figure out which player has won.
But we can use this concept of a game to construct a quite general
mathematical game $\Gamma_{p.i.}(f)$.
\begin{quote}
Let there be given two finite sets $X$ and $Y$, an integer $n$,
and a function $f$ assigning
to each sequence $w$ of length $n$ of pairs $(x_i,y_i) \in X\times Y$,
a {\em payoff} $f(w) \in \R$.
Two players are playing against each other.
Each player, in turn, makes a move by selecting an element $x_1 \in X$ or
$y_1 \in Y$, respectively, and announcing his or her selection.
Then they each in turn make a second move, and a third move, and continue
making moves until $n$ rounds have been played.
This generates a sequence $w$ of length $n$ of pairs $(x_i, y_i) \in X\times Y$.
Then they stop, and player II pays player I the amount $f(w)$.
\end{quote}
With the right choices for $X$, $Y$, $n$ and $f$, the game $\Gamma_{p.i.}(f)$
can `emulate' the game of Chess.
For if we let $X$ and $Y$ be the set of all possible chess moves,
and $n=6350$\footnote{The fifty-move rule is a rule in chess stating that
if no piece has been captured and no pawn has been moved for fifty turns,
the game is a draw.
Under the fifty-move rule, a game of chess can last a maximum of 6350 moves.},
then a sequence $w$ corresponds to a finished game of chess.
We now set $f(w)=1$, $f(w)=0$, or $f(w)=\frac{1}{2}$, depending on whether
the corresponding game is a win for White, a win for Black, or a draw.\footnote{
If $w$ does not correspond to a {\em legal} chess game, we count it
as a win for White if the first illegal move is made by Black, and vice versa.}
And voil\`a, we have our Chess emulator.

But Chess is not the only game that can be `emulated' in this manner.
The same can be done with Noughts-and-Crosses, Connect-Four, Go and Checkers.
In general, the games $\Gamma_{p.i.}(f)$ can emulate
any game $G$ that has the following properties:
\begin{itemize}
\item There are two players.
\item There is no element of chance
\item Moves are essentially made by selecting them and announcing them.
\item There is no hidden information:
a player knows all the moves made so far when making her current move,
and there is nothing going on simultaneously either ({\em Perfect Information}).
\item If one player loses (a certain amount) the other player wins
(that same amount) ({\em Zero-Sum}).
\item The game can last no more than a certain number of rounds 
({\em Finite Duration}).
\item There is a maximum number of alternatives each player can select from
({\em Finite Choice-of-Moves}).
\end{itemize}
Any results for the games $\Gamma_{p.i.}(f)$ apply to
all the games with these properties. \\
\\
David Blackwell described the concept of a {\em strategy} as \cite{Blackwell-3}:
\begin{quote}
Imagine that you are to play the White pieces in a single game of chess,
and that you discover you are unable to be present for the occasion.
There is available a deputy, who will represent you on the occasion,
and who will carry out your instructions exactly,
but who is absolutely unable to make any decisions of his own volition.
Thus, in order to guarantee
that your deputy will be able to conduct the White pieces throughout the game,
your instructions to him must envisage every possible circumstance
in which he may be required to move,
and must specify, for each such circumstance, what his choice is to be.
Any such complete set of instructions constitutes what we shall call a
{\em strategy}.
\end{quote}
Thus, a strategy for a given player in a given game consists of a specification,
for each position in which he or she is required to make a move,
of the particular move to make in that position.
In turn, a position can be specified by the moves made to get to that position.
If we apply this to the game $\Gamma_{p.i.}(f)$, a strategy becomes a function
 from the set of sequences of length $< n$ of pairs $(x_i, y_i) \in X\times Y$,
to the set of possible selections $X$, $Y$ respectively.

Given strategies for each of the players, the outcome of the game is determined:
each move follows from the current position and the strategy of the player whose
turn it is to move, and determines the next position.
So, the totality of all the decisions to be made can be described by a single
decision - the choice of a strategy.
This is the {\em normal form} of a game:
the two players independently make a single move,
which consists of selecting a strategy,
and then payoff is calculated and made.

Of course, there are good strategies and bad strategies.
The {\em value} of a strategy for a given player
is the result of that strategy against the best counterstrategy.
The {\em value} of a game for a given player is
the best result that that player can guarantee,
i.e. the value of that player's best strategy.
A game is called {\em determined} if its value is the same for both players.
That value is the result that will occur if both players are playing
perfectly.\footnote{In more general cases, we allow $\epsilon$-approximation,
i.e. a game is {\em determined} iff there exists a value $v$
such that for any $\epsilon>0$,
the two players have strategies guaranteeing them a payoff of
at least $v-\epsilon$ or at most $v+\epsilon$, respectively.}

Victor Allis \cite{Connect-Four} recently demonstrated that in a game
of Connect-Four, the first player can win, i.e. has a strategy that wins
against any counterstrategy.
And countless persons throughout the ages have independently discovered
that in the game of Noughts-and-Crosses, both players can force a draw.
These are both examples of determinacy.
It can be shown (using induction)
that any game $\Gamma(f)$ as defined above is determined,
and hence any game with all of the properties mentioned above is determined.
In the case of Go, Chess, and Checkers, this means that
either one of the players has a winning strategy,
or both players have a drawing strategy. \\
\\
Now consider the game of Scissors-Paper-Stone.
In this game, the two players simultaneously `throw' one of three symbols:
`Stone' (hand balled in a fist), `Paper' (hand flat with the palm down)
or `Scissors' (middle and forefinger spread, pointing forwards).
If both players throw the same symbols, the result is a draw;
otherwise, Paper beats Stone, Stone beats Scissors, and Scissors beats Paper
(the reason being that ``Paper wraps Stone, Stone blunts Scissors,
and Scissors cut Paper'').
In this game, the players do not make moves in turn, but simultaneously.
In other words, both players make moves, and neither player knows
what move the other is making.
This is an example of a game of {\em Imperfect Information}.

The strategy `Throw Stone' is a losing strategy,
because it loses against the counterstrategy `throw Paper'.
The same can be said for {\em any} strategy of the type `throw {\em this}',
for both players.
So in terms of the concept of strategy described above, this game is not
determined.
On the other hand, consider the `strategy'
 `throw Scissors $1/3$ of the time, throw Paper $1/3$ of the time,
 and throw Stone the remaining $1/3$ of the time'.
Against any other strategy, this strategy loses, draws and wins 
$1/3$ of the time each, for an `average result' of a draw.
This strategy does not fit in the concept of strategy given above,
but it is clearly worth considering.

Strategies of this new type are called {\em mixed} strategies, as opposed
to the old type of strategies, the {\em pure} strategies.
A mixed strategy for a given player in a given game consists of a specification,
for each position in which he or she is required to make a move,
of the {\em probability distribution} to be used to determine
what move to make in that position.\footnote{
Standard game theory defines a mixed strategy as a probability
distribution on pure strategies, but the above definition can be shown to be
equivalent to that one.}
Given mixed strategies for each of the players, 
the outcome of the game is not determined,
but we can calculate the {\em probability} of each outcome.
If we assign values to winning and losing 
(`the loser pays the winner one dollar'),
then we can calculate the average profit/loss one player can expect to make
 from the other, playing those strategies.

The {\em value} of a mixed strategy is therefore the
{\em expected average result} against the best counterstrategy.
And a game is called determined if, for some value $v$,
one of the players has a strategy with which she can always expect to make
(on average) at least $v$, no matter what the other plays,
while the other player has a strategy with which he can always expect to lose
(on average) at most $v$, no matter what the other plays.
As before, it can be shown (using induction and a theorem of Von Neumann)
that all finite two-person zero-sum games with Imperfect Information
(i.e. the games with the properties mentioned above,
except that players make moves at the same time instead of one after the other)
 are determined.

All the games mentioned so far are of {\em finite} duration.
Let, as before, $X$ and $Y$ be two finite sets, and
let $f$ be a function assigning
to each countably infinite sequence $w$ of pairs $(x_i,y_i) \in X\times Y$,
a payoff $f(w) \in \R$.\footnote{
We tacitly assume $f$ to be bounded, as otherwise things get ugly.}
We first consider games of infinite duration and {\em perfect} information:
\begin{quote}
Two players are playing against each other.
Each player, in turn, makes a move by selecting an element $x_1 \in X$ or
$y_1 \in Y$, respectively, and announcing his or her selection.
Then they each in turn make a second move, and a third move, and continue
making moves for a {\em countably infinite} number of rounds.
This generates an infinite sequence $w$ of pairs $(x_i, y_i) \in X\times Y$.
`Then' they stop, and player II pays player I the amount $f(w)$.
\end{quote}
The problem with infinite games, of course,
is that the outcome is only known after an infinite number of moves,
and thus it is impractical to play the game as it is.
But our concept of a strategy as a specification of
which move to make in each position,
is still valid in the case of games of infinite duration.
And given strategies for both players
we can construct the infinite sequence of moves that will be played
(or the probability distribution thereof), 
and apply the payoff function to obtain our (expectation of the) outcome.
Hence we can still play the game in a fashion, by using its {\em normal form}.

The concepts of values and determinateness carry over as well.
But it is no longer provable that all such games are determined.
For some payoff functions $f$, such as bounded Borel-measurable functions $f$,
it has been proven that the infinite game of perfect information
$\Gamma_{p.i.}(f)$ is determined.
But using the Axiom of Choice,
a nonmeasurable payoff function $f$ can be constructed such that
$\Gamma(f)$ is not determined \cite{AxiomOfChoice}.
The axiom AD, the axiom that all games $\Gamma(f)$ are determined,
is widely used as an alternative to AC \cite{AD-1,AD-2}.

A game of infinite duration and {\em imperfect} information is similar, except
that both players make their $n^{th}$ move at the same time.
These games are called Blackwell games,
named after David Blackwell, the first one
to describe and study these games \cite{Blackwell}.
For Blackwell games, it has been proven that $\Gamma(f)$ is determined
for the case that $f$ is the indicator function of an open or $G_\delta$ set.
In this article I prove determinacy of $\Gamma(f)$ for the case
that $f$ is the indicator function of a $G_{\delta\sigma}$ set.
But the general case of Borel-measurable functions is still open.
\clearpage
\Section{Definitions, Lemmas and Terminology}

\subsection{Games, Strategies and Values}
The definitions in this subsection are fairly standard,
and merely formalize the intuitive concepts from the introduction.
The lemmas are all basic properties of game-values.
For reasons of conciseness, no proofs are given in this section.
\\ \ \\
In order to define what a Blackwell game is, we first need some sets.
Let $X$ and $Y$ be two finite, nonempty sets, and put $Z = X \times Y$.
An \em (infinite) play \em is a countably infinite sequence $w$ of pairs 
$(x,y) \in Z$.
We write $W$ for the set of all plays, i.e. $W=Z^\N$ \\

\begin{defn}
Let $f: W \to \R$ be a bounded Borel (measurable) function
(i.e. a bounded function such that $f^{-1}[O]$ is a Borel set
for every open set $O \subseteq \R$).
The \em Blackwell game $\Gamma(f)$ \em with \em payoff function $f$ \em
is the two-person zero-sum infinite game of imperfect information
 played as follows: 
Player I selects an element $x_1 \in X$ (\em makes the move $x_1$\em)
and, simultaneously, player II selects an element $y_1 \in Y$. 
Then both players are told $z_1 = (x_1, y_1)$, and the game \em is at \em
or \em has reached \em position $(z_1)$.
Then player I selects $x_2 \in X$ and, simultaneously, 
player II selects $y_2 \in Y$.
Then both players are told $z_2 = (x_2, y_2)$, and the game is at
position $(z_1, z_2)$.
Then both players simultaneously selects $x_3 \in X$ and $y_3 \in Y$, etc.
Thus they produce a play $w=(z_1, z_2, \ldots)$.
Then player II pays player I the amount $f(w)$, ending the game.
\end{defn}
A \em position \em or \em finite play \em (\em of length $k$\em) is a finite
sequence $p$ (of length $k$) of pairs $(x,y) \in Z$.
We write $P$ for the set of all positions, i.e. $P=Z^{< \omega}$. \\ \ \\
Some notation and terminology that we are going to use: \\
$w$ usually denotes an infinite play, $p$ denotes a finite play or position. \\
$p_{\mid n}$, $w_{\mid n}$ denote the sequences consisting 
of the first $n$ moves made in $p$, $w$ respectively
(counting a pair $(x,y)$ as one move). \\
$p\cat p'$, $p\cat w$ denote the sequences consisting of
the finite sequence $p$ followed by the finite sequence $p'$
or the infinite sequence $w$, respectively. \\
$\len(p)$ denotes the length of a finite sequence $p$. \\
$e$ denotes the position of length $0$, i.e. the empty sequence. \\
$W_n$ denotes the set of all finite plays of length $n$,
i.e. $W_n=Z^n$, for $n \in \N$. \\
$p \subset w$ denotes that $w_{\mid\mbox{\scriptsize len}(p)}=p$,
and we say that $w$ \em hits \em or \em passes through \em $p$. \\
$p \subset p'$ denotes that $p'_{\mid\mbox{\scriptsize len}(p)}=p$
and $p' \neq p$, and we say that
$p'$ \em follows \em $p$, and $p$ \em precedes \em $p'$. \\
$p \subseteq p'$ denotes that $p'_{\mid\mbox{\scriptsize len}(p)}=p$,
and we say that $p'$ \em follows or is equal to \em $p$. \\
$[p]$ denotes the set $\{ w \in W \mid w \supset p \}$
of all plays hitting the position $p$.\\
$[H]$ denotes the set $\{ w \in W \mid \exists p \in H : w \supset p \}$
of all plays hitting any position in a set of positions $H$.\\
We will sometimes write a sequence $((x_1, y_1), (x_2, y_2), \ldots)$
as $(x_1, y_1, x_2, y_2, \ldots)$.
$\Gamma(S)$ denotes the game $\Gamma(I_S)$, where $I_S$ is the indicator
function of $S \subseteq W$.

\begin{rem}
We give $W$ the usual topology by
letting the basic open sets be the sets of the form $[H]$
for some some set $H \subseteq W_n$ of positions of fixed length $n$.
Then the open subsets of $W$ are exactly those of the form $[H]$ 
for some set $H$ of positions. 
The $G_\delta$ subsets of $W$ are exactly those of the form
$ \{ w \in W \mid \#\{ p \in H \mid w\mbox{ hits }p \} = \infty \} $
for some set $H$ of positions.
Note that under this topology, $W$ is a compact space.
\end{rem}

\begin{defn}
A \em strategy \em for player I in a Blackwell game $\Gamma(f)$
is a function $\sigma$ assigning to each 
position $p$ a probability distribution on $X$.
More formally, $\sigma$ is a function $P \to [0, 1]^X$ satisfying 
$\forall p \in P: \sum_{x \in X} \sigma(p)(x) = 1$.

Analogously, a \em strategy \em for player II is a function $\tau$ assigning
to each position $p$ a probability distribution on $Y$.
\end{defn}

\begin{defn}
Let $\sigma$ and $\tau$ be strategies for players I, II in a Blackwell game
$\Gamma(f)$. $\sigma$ and $\tau$ determine a probability measure 
$\mu_{\sigma, \tau}$ on $W$, induced by 
\begin{equation}
\mu_{\sigma, \tau}[p] 
= P\{ w \mid w\mbox{ hits }p\}
=\prod_{i=1}^n \left(
\sigma(p_{\mid(i-1)})(x_i) \bullet \tau(p_{\mid(i-1)})(y_i) \right)
\end{equation}
for any position $p = (x_1, y_1, \ldots, x_n, y_n) \in P$.

The \em expected income \em of player I in $\Gamma(f)$,
if she plays according to $\sigma$ 
and player II plays according to $\tau$, is the expectation of $f(w)$ 
under this probability measure:
\begin{equation}
\expectst{f} = \int f(w)d\mu_{\sigma, \tau}(w)
\end{equation}
\end{defn}

\begin{defn}
Let $\Gamma(f)$ be a Blackwell game.
The \em value \em of a strategy $\sigma$ for player I in $\Gamma(f)$ is
the expected income player I can guarantee if she plays according to $\sigma$.
Similarly, the \em value \em of a strategy $\tau$ for player II in
$\Gamma(f)$ is the amount to which player II can restrict player I's income
if he plays according to $\tau$.
I.e.
\begin{eqnarray}
\val(\sigma\mbox{ in }\Gamma(f)) & = & \inf_{\tau}\expectst{f} \\
\val(\tau\mbox{ in }\Gamma(f)) & = & \sup_{\sigma}\expectst{f}
\end{eqnarray}
\end{defn}

\begin{defn}
Let $\Gamma(f)$ be a Blackwell game. The \em lower value \em of $\Gamma(f)$
is the smallest upper bound on the income that player I can guarantee.
Similarly, the \em upper value \em of $\Gamma(f)$ is the largest lower bound
on the restrictions player II can put on player I's income.
I.e.
\begin{eqnarray}
\lowerval(\Gamma(f)) =   \sup_{\sigma}\val(\sigma\mbox{ in }\Gamma(f))
		   & = & \sup_{\sigma}\inf_{\tau}\expectst{f} \\
\upperval(\Gamma(f)) =   \inf_{\tau}\val(\tau\mbox{ in }\Gamma(f))
		   & = & \inf_{\tau}\sup_{\sigma}\expectst{f}
\end{eqnarray}
Clearly, for all games $\Gamma(f)$,
$\lowerval(\Gamma(f)) \leq \upperval(\Gamma(f))$.
If $\upperval(\Gamma(f)) = \lowerval(\Gamma(f))$, then $\Gamma(f)$ is called
\em determined\em, and we write $\val(\Gamma(f)) = 
\upperval(\Gamma(f)) = \lowerval(\Gamma(f))$.
\end{defn}

\begin{defn}
Let $\Gamma(f)$ be a Blackwell game, and let $\epsilon>0$.
A strategy $\sigma$ for player I in $\Gamma(f)$ is 
\em optimal \em if 
$\val(\sigma\mbox{ in }\Gamma(f)) = \lowerval(\Gamma(f))$.
A strategy $\sigma$ for player I in $\Gamma(f)$ is 
\em $\epsilon$-optimal \em if 
$\val(\sigma\mbox{ in }\Gamma(f))>\lowerval(\Gamma(f))-\epsilon$.
Similarly, a strategy $\tau$ for player II in $\Gamma(f)$ is
\em optimal \em if 
$\val(\tau\mbox{ in }\Gamma(f)) = \upperval(\Gamma(f))$, and
\em $\epsilon$-optimal \em if 
$\val(\tau\mbox{ in }\Gamma(f)) < \upperval(\Gamma(f))+\epsilon$.
\end{defn}
Some basic properties of these values are:

\begin{lem}
\label{Payoff-Compare-Lemma}
Let $f$, $g$ be two payoff functions such that for all $w \in W$,
$f(w) \leq g(w)$.
Then 
$\lowerval(\Gamma(f)) \leq \lowerval(\Gamma(g))$ and
$\upperval(\Gamma(f)) \leq \upperval(\Gamma(g))$.
\end{lem}

\begin{lem}
\label{Payoff-Linear-Lemma}
Let $f$ be a payoff function, and let $a, c \in \R, a \geq 0$.Then
\linebreak[4]
$\lowerval(\Gamma(af+c)) = a\,\lowerval(\Gamma(f)) + c$ and
$\upperval(\Gamma(af+c)) = a\,\upperval(\Gamma(f)) + c$.
\end{lem}

\begin{lem}
\label{Payoff-Switched-Lemma}
Let $f$ be a payoff function, and let $\switched{f}: (Y \times X)^\N \to \R$
be defined by
\begin{equation}
\switched{f}((y_1, x_1), (y_2, x_2), \dots) = f((x_1, y_1), (x_2, y_2), \dots)
\end{equation}
Then 
\begin{eqnarray}
\lowerval(\Gamma(-f)) & = & -\upperval(\switched{\Gamma}(\switched{f})) \\
\upperval(\Gamma(-f)) & = & -\lowerval(\switched{\Gamma}(\switched{f}))
\end{eqnarray}
where $\switched{\Gamma}(\switched{f})$ is the Blackwell game with payoff
function $\switched{f}$ in which player I selects moves from $Y$ and
player II selects moves from $X$.
\end{lem}

\begin{lem}
\label{Limit-Pointwise-Lemma}
Let $(f_i)_i$ be a sequence of functions $f_i: W \to [a, b]$ such that 
$(f_i)_i$ converges pointwise to a function $f: W \to [a, b]$.
Then for any two strategies $\sigma$, $\tau$, 
$\lim_{i \to \infty} \expectst{f_i} = \expectst{f}$
\end{lem}

\subsection{Starting and Stopping}

\begin{defn}
Let $f: W \to \R$ be a bounded Borel function, and 
$p=((x_1, y_1), (x_2, y_2), \ldots, (x_n, y_n))$ a position.
The \em subgame $\Gamma(f, p)$ starting from position $p$ \em is
the game played like $\Gamma(f)$,
except that the players start at round $n+1$, and the first $n$ moves are
supposed to have been $x_1, y_1, x_2, y_2, \ldots, x_n, y_n$.
The game $\Gamma(f, p)$ is played exactly the same as the game $\Gamma(g)$,
where $g$ is the payoff function defined by $g(w)=f(p\cat w)$.

As before, strategies $\sigma$ and $\tau$ determine a probability measure 
$\mu_{\sigma, \tau\mbox{ \scriptsize in }\Gamma(f,p)}$ on $W$.
This measure is equal to the conditional probability measure obtained from
$\mu_{\sigma, \tau}$ given $[p]$, i.e.
\begin{equation}
\mu_{\sigma, \tau\mbox{ \scriptsize in }\Gamma(f,p)}(S) =
             \frac{\mu_{\sigma, \tau}(S \cap [p])} {\mu_{\sigma, \tau}[p]}
\mbox{ \ \ if } \mu_{\sigma, \tau}[p] = 0
\end{equation}
The expected income of player I, the value of a strategy $\sigma$, etc. are
defined for the games $\Gamma(f,p)$ in the same manner as for the games 
$\Gamma(f)$.
 \end{defn}

\begin{defn}
A \em stopping position \em in a Blackwell game $\Gamma(f)$ is a position
$p$, such that for all plays $w, w' \in [p]$, $f(w)=f(w')$.
We will denote this value by $f(p)$.
A \em stopset \em in a Blackwell game $\Gamma(f)$ is a set $H$ of stopping
positions, such that no stopping position $p \in H$ precedes another
stopping position $p' \in H$.
\end{defn}
We will often define a payoff function $f$ using the following format:
\begin{eqnarray*}
\payoffdef{f}{H}{\mbox{formula1}}{\mbox{formula2}}
\end{eqnarray*}
where $H$ is a set of positions such that no position
$p \in H$ precedes another position $p' \in H$.
Then H is a stopset in the game $\Gamma_H(f)$.
\begin{rem}
\label{Undefined-Strategy-Remark}
If $p$ is a stopping position, any moves made at or after $p$ will not affect
the outcome of the game. It is often convenient to assume that both players
will stop playing if a stopping position is reached. If $\Gamma(f)$ is a
Blackwell game, and $H$ is a stopset, we write $\Gamma_H(f)$ to
explicitly denote that players stop playing at the positions in $H$.
In this case, we only require strategies to be defined on nonstopping positions.
Similarly, with respect to a subgame $\Gamma(f, p)$, we only require
strategies to be defined on positions that are following or equal to $p$.
In fact it is occasionally necessary to assume that a strategy is \em not \em
defined on positions outside the subgame proper,
for instance to combine strategies for different subgames into one big strategy.
\end{rem}
Using stopsets, a finite game can be treated as a special type of infinite game.
For finite games, we have determinacy, as well as a kind of continuity of
the value function.
\begin{defn}
Let $\Gamma(f)$ be a Blackwell game.
If, for some $n$, all positions in $W_n$ are stopping positions, then
$\Gamma(f)$ is called \em finite (of length $n$)\em.
If $\Gamma(f)$ is finite, we can stop after playing $n$ rounds, and we will
denote this by writing $\Gamma_n(f)$.
\end{defn}

\begin{thm}[Von Neumann's Minimax Theorem\cite{Neumann}]
\label{One-Move}
Let $\Gamma_1(f)$ be a finite one-round Blackwell game (i.e. of length 1).
Then $\Gamma_1(f)$ is determined, and both players have optimal strategies.
\end{thm}

\begin{thm}
\label{Finite-Theorem}
Let $\Gamma_n(f)$ be a finite Blackwell game of length $n$.
Then $\Gamma_n(f)$ is determined, and both players have optimal strategies.
\end{thm}

\begin{lem}
\label{Limit-Finite-Lemma}
Let $n \in \N$.
Let $(f_i)_i$ be a sequence of payoff functions $f_i: W_n \to [a, b]$
such that $(f_i)_i$ converges to a payoff function
$f: W_n \to [a, b]$.
Then $ \val(\Gamma_n(f)) = \lim_{i \to \infty} \val(\Gamma_n(f_i)) $.
\end{lem}

\subsection{Equivalent Truncated Subgames}

In games like Chess, Go, or even Risk or Monopoly,
a player is usually allowed to give up if he has no hope of winning.
He doesn't have to play it out
in the hope that the other player will make a mistake.
Two players can agree beforehand to stop in certain positions, and pay
out the value of the game at that position rather than continue playing.
Provided their assessment of that value is accurate, this does
not change the value of the total game.
We will call a game resulting from such an alteration 
a {\em truncated subgame}.

\begin{defn}
\label{Equivalent-Definition}
Let $f$, $g$ be two payoff functions, and $H$ a stopset in $\Gamma(g)$.
$\Gamma_H(g)$ is an \em equivalent truncated subgame \em of $\Gamma(f)$
(\em truncated at $H$\em),
if for any play $w \not\in [H]$, $f(w)=g(w)$,
and for any $p \in H$, $g(p) = \val(\Gamma(f, p))$.
$\Gamma_H(g)$ is a \em truncated subgame,
equivalent for player I [player II]\em,
if for any play $w \not\in H$, $f(w)=g(w)$,
and for any $p \in H$, $g(p) = \lowerval(\Gamma(f, p))$
[$g(p) = \upperval(\Gamma(f, p))$].
In all three cases,
$\Gamma(f)$ is called an \em extension \em of $\Gamma_H(g)$.
\end{defn}
Note that $\Gamma_H(g)$ is an equivalent truncated subgame of $\Gamma(f)$ iff
it is a truncated subgame equivalent for both player I and player II.

\begin{lem}
\label{Uppervalue-Lemma}
\label{Lowervalue-Lemma}
Let $\Gamma(f)$ be a Blackwell game,
and let $\Gamma_H(g)$ be a truncated subgame of $\Gamma(f)$,
truncated at a set of positions $H$, equivalent for player I [for player II].
Then $\lowerval(\Gamma(f))=\lowerval(\Gamma_H(g))$ 
[$\upperval(\Gamma(f))=\upperval(\Gamma_H(g))$]. Furthermore,
for any $\epsilon>0$, 
any $\epsilon$-optimal strategy for player I [player II] in $\Gamma_H(g)$
(if it is undefined on all positions at or after positions in $H$) 
can be extended to an $\epsilon$-optimal strategy for player I [player II] 
in $\Gamma(f)$.
\end{lem}
\prc
We find an $\epsilon$-optimal strategy for the truncated subgame $\Gamma_H(g)$,
and for the appropriate $\delta$,
$\delta$-optimal strategies for the games $\Gamma(f, p)$ starting at 
positions $p \in H$, i.e. the positions where $\Gamma_H(g)$ stops.
Then we tie them together, and calculate how well the combination strategy
performs against opposing strategies.

\begin{cor}
\label{Value-Corollary}
Let $\Gamma(f)$ be a Blackwell game,
and let $\Gamma_H(g)$ be an equivalent truncated subgame of $\Gamma(f)$
(truncated at $H$).
If $\Gamma_H(g)$ is determined, then $\Gamma(f)$ is determined, and
$\val(\Gamma(f))=\val(\Gamma_H(g))$. Furthermore, for any $\epsilon>0$,
any $\epsilon$-optimal strategy for player I or player II in $\Gamma_H(g)$ 
can be extended to an $\epsilon$-optimal strategy for player I or player II 
in $\Gamma(f)$.
\end{cor}

\begin{cor}
\label{Lowervalue-Compare-Corollary}
\label{Uppervalue-Compare-Corollary}
\label{Value-Compare-Corollary}
Let $\Gamma(f), \Gamma_H(g)$ be Blackwell games. 
If for any $p \in H$,
$g(p) \leq \lowerval(\Gamma(f,p))$,
and for any $w \not\in [H]$, $g(w) \leq f(w)$,
then $\lowerval(\Gamma_H(g)) \leq \lowerval(\Gamma(f))$.
Similarly for the value and the upper value, and for $\geq$ instead of $\leq$.
\end{cor}
Truncated subgames may be nested.
If we have a nested series of truncated subgames, then we may extend
a strategy for the smallest subgame to a strategy for {\em all} subgames.
This allows us to approximate complicated games with a series of simpler,
truncated subgames, obtain a strategy that is ($\epsilon$-)optimal in all
the subgames.
The final lemma in this section allows us to prove results for that strategy
in the original game.

\begin{defn}
Let, for $n \in \N$, $f_n$ be a payoff function, and
$H_n$ a set of stopping positions in $\Gamma(f_n)$.
If for all $n \in \N$, $\Gamma_{H_n}(f_n)$ is a truncated subgame of 
$\Gamma_{H_{n+1}}(f_{n+1})$, 
and equivalent to $\Gamma_{H_{n+1}}(f_{n+1})$ [for player I, II],
then the series of games $(\Gamma_{H_n}(f_n))_{n \in \N}$ is
called a \em nested series of equivalent truncated subgames \em
[equivalent for player I, II].
\end{defn}

\begin{lem}
\label{Nested-Lowervalue-Lemma}
\label{Nested-Uppervalue-Lemma}
Let $(\Gamma_{H_i}(g_i))_{i \in \N}$ be a nested series of truncated games
equivalent for player I [player II]. 
Then all the games have the same lower value [upper value].
Furthermore, we can find a strategy for player I [player II] 
that is $\epsilon$-optimal in all the games $\Gamma_{H_i}(g_i)$.
\end{lem}
\prc
Basically, we apply Lemma \ref{Lowervalue-Lemma} a number of times
and use induction.
The proof is straightforward,
except for a slight complication involving the domain on which the strategies
are defined.
This complication is solved using the observations
that if we truncate a game, any stopping position remains a stopping position,
and that strategies can be assumed to be undefined on stopping positions.

\begin{cor}
\label{Nested-Value-Corollary}
Let $(\Gamma_{H_i}(g_i))_{i \in \N}$ be a nested series of equivalent truncated
subgames. If $\Gamma_{H_0}(g_0)$ is determined, then all the
games are determined, and all the games have the same value.
Furthermore, we can find strategies for player I and player II 
that are $\epsilon$-optimal in all the games $\Gamma_{H_i}(g_i)$.
\end{cor}

\begin{rem}
\label{Strategy-Extension-Optimal-Remark}
If the component games involved all have optimal strategies,
then we can extend optimal strategies with optimal strategies
to optimal strategies,
i.e. drop the $\epsilon$ in the above lemmas and corollaries.
\end{rem}

\section{Determinateness Results}

\subsection{Generalized Open Games}

In this subsection we prove determinacy of 
a class of `generalized open games', where
payoff for a play is calculated as the supremum of values associated
with the positions hit in the play.
In addition we derive a result for these and open games
comparable to the compactness of $W$.

\begin{thm}
\label{Open-Function-Theorem}
Let $u: P \to \R$ be a bounded function, and let $f: W \to \R$ be 
the payoff function defined by $f(w) = \sup_{j \in \N} u(w_{\mid j})$.
Then $\Gamma(f)$ is determined, and
\begin{equation} \val(\Gamma(f)) = \lim_{n \to \infty} \val(\Gamma_n(f_n))
\end{equation}
where $f_n(w) = \sup_{j \leq n} u(w_{\mid j})$.
\end{thm}
\prc
Showing that $\lim_{n \to \infty} \val(\Gamma_n(f_n)) $ exists and is not 
greater than the lower value of $\Gamma(f)$ is not difficult.
To show that it is not less than the upper value,
we approximate $\Gamma(f)$ with a collection of finite
auxiliary games $\Gamma_n(g_n)$ such that the payoff at the stopping 
positions is an \em estimate \em of the value of the game at that point. 
We then show that these auxiliary games form a 
nested series of equivalent finite truncated subgames.
This allows us to find a strategy that is optimal in
each of the truncated subgames. This strategy is also a strategy
in the game $\Gamma(f)$, and has a value in $\Gamma(f)$ equal to 
 $\lim_{n \to \infty} \val(\Gamma_n(f_n))$.
\\ \ \\
\proof
Without loss of generality we may assume that the function 
$u$ has range $[0, 1]$.
For any $p \in P$, and any $n \in \N$, the game $\Gamma_n(f_n, p)$ is finite
(of length $\leq n$), and thus determined.
It is easily seen that $f_0 \leq f_1 \leq f_2 \leq \ldots \leq f \leq 1$.
Consequently, for any $p \in P$,
\begin{equation}
\val(\Gamma_0(f_0, p)) \leq \val(\Gamma_1(f_1, p)) \leq \val(\Gamma_2(f_2, p))
\leq \ldots \leq \lowerval(\Gamma(f, p)) \leq 1
\end{equation}
For all $p \in P$, $\lim_{k \to \infty} \val(\Gamma_k(f_k, p))$ exists,
since all monotone non-decreasing bounded sequences converge.
Furthermore, for all $p \in P$,
\begin{equation} 
\lim_{n \to \infty} \val(\Gamma_n(f_n, p)) \leq \lowerval(\Gamma(f, p))
\end{equation}
Define for any $n \in \N$ the payoff function $g_n: W_n \to [0, 1]$ by
\begin{equation}
g_n(p) = \lim_{k \to \infty} \val(\Gamma_k(f_k, p))\mbox{ for }p \in W_n
\end{equation}
Then for all $p \in W_n$, $g_n(p) \geq \val(\Gamma_n(f_n, p)) = f_n(p)$.

Furthermore, the games $\Gamma_n(g_n)$ form a nested series of equivalent
truncated subgames.
For fix $n \in \N$, $p \in W_n$. 
Define for $k \in \N$, $h_{n+1, k}: W_{n+1} \to \R$ by
$h_{n+1, k}(p') = \val(\Gamma_k(f_k, p'))$ for $p' \in W_{n+1}$.
Then
\begin{eqnarray}
g_n(p) & = & \lim_{k \to \infty} \val(\Gamma_k(f_k, p)) \\
\label{equation-h}
       & = & \lim_{k \to \infty} \val(\Gamma_{n+1}(h_{n+1, k}, p)) \\
\label{equation-limit-finite}
       & = & \val(\Gamma_{n+1}(\lim_{k \to \infty} h_{n+1, k}, p)) \\
       & = & \val(\Gamma_{n+1}(g_{n+1}, p))
\end{eqnarray}
(equation (\ref{equation-h}) follows from Corollary \ref{Value-Corollary}, and
equation (\ref{equation-limit-finite}) follows from
Lemma \ref{Limit-Finite-Lemma}
as $W_{n+1}$ is finite).

Since $(\Gamma_n(g_n))_{n \in \N}$ is a nested series of equivalent truncated
subgames, by Corollary \ref{Nested-Value-Corollary}
the games $\Gamma_n(g_n)$ all 
have the same value, say $v$.
Also, we can find a strategy for player II 
that is $\epsilon$-optimal in all the games $\Gamma_n(g_n)$,
and since all the games $\Gamma_n(g_n)$ are finite and hence have optimal strategies,
by Remark \ref{Strategy-Extension-Optimal-Remark}
we can even find a strategy that is optimal in all the games $\Gamma_n(g_n)$.
So let $\tau$ be such a strategy.
Then for any strategy $\sigma$, and any $n \in \N$,
\begin{equation} 
\expectfin{n}{g_n} \leq \val(\Gamma_n(g_n)) = v
\end{equation}
Now let $\sigma$ be any strategy for player I in $\Gamma(f)$.
Then
\begin{eqnarray}
\lefteqn{ \expectst{f} } \nonumber \\
& = & \lim_{n \to \infty} \expectfin{n}{f_n} \\
& \leq & \lim_{n \to \infty} \expectfin{n}{g_n} \\
& \leq & v
\end{eqnarray}
So 
\begin{equation}
\upperval(\Gamma(f)) \leq \val(\tau \mbox{ for player II in } \Gamma(f))
\leq v
\end{equation}
But also
\begin{equation}
v = \val(\Gamma_0(g_0)) = g_o(e) = \lim_{k \to \infty} \val(\Gamma_k(f_k))
\leq \lowerval(\Gamma(f))
\end{equation}
Therefore,
\begin{equation}
\upperval(\Gamma(f)) = \lowerval(\Gamma(f)) = \lim_{k \to \infty}
\val(\Gamma_k(f_k))
\end{equation}
\eop

\begin{cor}
\label{Open-Set-Corollary}
Let $O$ be an open set.
Then $\Gamma(O)$ is determined.
\end{cor}
\proof
There exists a set of positions $H$ such that $O=[H]$.
Then for all $w \in W$, $I_O(w) = \sup_{n \in \N} I_H(w_{\mid n})$.
Applying Theorem \ref{Open-Function-Theorem} yields the corollary.
\eop

\begin{cor}
\label{Open-Union-Corollary}
Let $O = \bigcup_i O_i$ be the union of open sets.
Then $\val(\Gamma(O)) = \lim_{n \to \infty} \val(\Gamma(\bigcup_{i \leq n} O_i))$.
\end{cor}
\proof
As the union of open sets, $O$ is open,
and hence there is a set of positions $H$ such that $O=[H]$, i.e.
\begin{equation}
O = \{ w \in W \mid \exists p \in H : p \subset w \}
\end{equation}
Define the basic open sets $B_j \subseteq O$ by
\begin{equation}
B_n = \{ w \in W \mid \exists p \in H : p \subset w \wedge \mbox{len}(p) \leq n \}
\end{equation}
then for all $w \in W$, 
\begin{eqnarray}
I_O(w) & = & \sup_{j \in \N} I_H(w_{\mid j}) \\
I_{B_n}(w) & = & \sup_{j \leq n} I_H(w_{\mid j})
\end{eqnarray}
so applying Theorem \ref{Open-Function-Theorem}, we find that
\begin{equation}
\val(\Gamma(O)) = \lim_{n\to \infty} \val(\Gamma_n(B_n))
\end{equation}
For each $m \in \N$, $B_m$ is a closed set covered by the open
sets $(O_i)_{i \in \N}$.
So by the compactness of $W$ there is for each $m \in \N$ a
$n_m \in \N$ such that $B_m \subseteq \bigcup_{i=1}^{n_m} O_i$. 
Then for all $n \geq n_m$, 
\begin{equation}
\val(\Gamma_m(B_m)) \leq \val(\Gamma(\bigcup_{i=1}^n O_i)) \leq \val(\Gamma(O)) 
\end{equation}
The corollary follows immediately. 
\eop

\begin{cor}
Let $f$ be a continuous function.
Then $\Gamma(f)$ is determined.
\end{cor}
\proof
As $W$ is compact, and $f$ is continuous, $f[W]$ is compact, and hence bounded.
Define $u: P \to \R$ by $u(p) := \inf_{w\in[p]} f(w)$.
Then $u$ is well-defined and bounded, and by the continuity of $f$,
$f(w) = \sup_{n \in \N} u(w_{\mid n})$ for all $w \in W$.
Applying Theorem \ref{Open-Function-Theorem} yields the corollary.
\eop

\begin{rem}
\label{Open-Local-Strategy-Remark}
In the case of open games (or generalized open games, described later)
there is an optimal strategy for player II.
This strategy can be described as `at every position player II plays
the optimal one-round strategy, looking at the values the game has
\em for player I \em from all positions directly following that one'.
However, for player I there does not always exist an optimal strategy,
as the following example shows.
\end{rem}

\begin{ex}
\label{Anti-Local-Strategy-Example}
Consider the following Blackwell game.
Each round, both players say either `Stop' or `Continue'.
If both players say `Continue', then play continues.
Otherwise, the game halts:
player II wins (payoff $0$) if both players said `Stop',
while player I wins (payoff $1$) if only one of the players said `Stop'.
If play continues indefinitely, and neither player ever says `Stop',
then payoff is $0$, i.e. player II wins.

This is clearly an open game.
An interpretation of this game is, that player II tries to guess
on which round player I will say `Stop', and tries to match her.
If player II guesses wrong, i.e. says `Stop' too soon or not soon enough,
then player I wins, if player II guesses right, then he wins.

A strategy of value $1-\frac{1}{n}$ for player I is, to select at random a
number $i$ between $1$ and $n$, and say `Stop' on round $i$.
Translated to the standard format for strategies, this becomes:

on round $1$, say `Stop' $\frac{1}{n}$ of the time,

on round $2$, if not yet stopped, say `Stop' $\frac{1}{n-1}$ of the time,

on round $3$, if not yet stopped, say `Stop' $\frac{1}{n-2}$ of the time,

\hspace{.5in} \vdots

on round $n$, if not yet stopped, say `Stop' $\frac{1}{1}$ of the time.

Hence, the value of this game is $1$. In fact, the value of this game at
any position in which game has not yet ended is $1$.
But there exists no optimal strategy of value $1$.
For suppose there exists such a strategy, of value $1$.
Then on any round (in which play has not yet ended),
the chance that player I will say `Stop' in that round is $0\%$.
For otherwise, the strategy would not score $100\%$ against the
counterstrategy that player II says `Stop' on that round.
But then, player I will never say `Stop', and this strategy will lose
against the counterstrategy that player II never says `Stop'.
So any strategy for player I has value strictly less than $1$,
although there are strategies with values arbitrarily close to $1$.
This game is an example of
a game in which one of the players has no optimal strategy.
\end{ex}

\subsection{$G_\delta$-sets}

Davis' proof of determinacy for $G_{\delta\sigma}$ games of perfect information
\cite{Davis} is based upon the idea of `imposing restrictions' on
the range of moves player II can make.
I.e. certain moves are declared `forbidden', or a loss for player II,
in such a way that 
(a) if player I did not have a win before, she does not get a win now, and
(b) a particular $G_\delta$ set is now certain to be avoided.
By applying this to all the $G_\delta$ subsets of a $G_{\delta\sigma}$ set,
and using compactness, he shows that
if player I cannot force the resulting sequence to be in one of the
$G_\delta$ sets, player II can force the resulting sequence to be outside
all of them.

The union of all the sequences in which one of the `forbidden' moves is played, is an open set that contains the $G_\delta$ set in question.
One way of looking at Davis' proof is,
that he enlarges each of the $G_\delta$ sets to an open set
without increasing the (lower) value of the game,
in order to be able to apply determinacy of open games.
 
In this subsection, we show that this holds (in a fashion) for Blackwell games,
i.e. that a $G_\delta$ set can be `enlarged' to an open set
without increasing the lower value of the game 
by more than an arbitrarily small amount,
even in the presence of a `background function',
a payoff function for those sequences that are not in the $G_\delta$ set.

\begin{thm}
\label{G-Delta-And-Function-Theorem}
Let $f: W \to [0,1]$ be a measurable function and let $D$ be a $G_\delta$ set.
Then
\begin{equation}
\lowerval(\Gamma(\max(f, I_D))) =
  \inf_{O \supseteq D, O\mbox{ \scriptsize open}}
    \lowerval(\Gamma(\max(f, I_O)))
\end{equation}
\end{thm}

\prc
We define a collection of
auxiliary games $\Gamma_{H_i}(g_i)$ of the game $\Gamma(\max(f, I_D))$, in
which the amount player I gets at a stopping position $p$ is
an \em estimate \em for the value of $\Gamma(\max(f, I_D))$ at position $p$,
namely $\inf_{O \supseteq D, O\mbox{ \scriptsize open}}
\val(\Gamma(\max(f, I_O), p))$.
We then show that these auxiliary games form a 
nested series of finite truncated subgames, equivalent for player I.
This allows us to find a strategy that is $\epsilon$-optimal in
each of the truncated subgames. This strategy is also a strategy
in the game $\Gamma(\max(f,I_D))$, and has the required value, proving one
side of the equation. The other side is trivial.
\\ \ \\
\proof
Put $v=\inf_{O \supseteq D, O\mbox{ \scriptsize open}} \lowerval(\Gamma(\max(f, I_O)))$.
For any $G_\delta$ set $D$ we can find a set of positions $H$, such
that $ D = \{w \in W \mid \#\{ p \in H \mid p \subset w \} = \infty \} $.
We may assume that $e \in H$.
\\ \ \\
Define for any $i \in \N$,
\begin{equation} H_i := \{ p \in H \mid
\mbox{there are exactly $i$ positions $p'$ in $H$ strictly preceding $p$} \}
\end{equation}
Define for any $i \in \N$ the payoff functions $g_i$, $h_i$ by
\begin{eqnarray}
\payoffdef{g_i}{H_i}
{\inf_{O \supseteq D, O\mbox{ \scriptsize open}} \lowerval(\Gamma(\max(f, I_O), p))}{f(w)} \\
\payoffdef{h_i}{H_i}{1}{f(w)}
\end{eqnarray}
First, the games $\Gamma_{H_i}(g_i)$ form a nested series of 
truncated subgames equivalent for player I.
For let $i \in \N$, and fix $p \in H_i$.
Let $O \supseteq D$, then for any $p' \in H_{i+1}$ such that $p' \supseteq p$,
$\lowerval(\Gamma(\max(f,I_O), p')) \geq g_{i+1}(p')$,
and for any $w \supset p$ such that $w \not\in [H_{i+1}]$,
$\max(f,I_O)(w) \geq f(w)=g_{i+1}(w)$.
Hence by Corollary \ref{Lowervalue-Compare-Corollary},
for any $O \supseteq D$, 
$\lowerval(\Gamma(\max(f,I_O), p)) \geq \lowerval(\Gamma_{H_{i+1}}(g_{i+1}, p)) $.
Therefore,
\begin{equation}
g_i(p) \geq \lowerval(\Gamma_{H_{i+1}}(g_{i+1}, p))
\end{equation}
On the other hand,
for any $\epsilon>0$ we can find, for each $p' \in H_{i+1}$, 
an open set $O_{p'} \supseteq D$ such that 
\begin{equation} 
\lowerval(\Gamma(\max(f,I_{O_{p'}}), p')) \leq g_{i+1}(p') + \epsilon
\end{equation}
Set $O = \bigcup_{p' \in H_{i+1}}([p']\cap O_{p'})$.
Then for all $p' \in H_{i+1}$,
$\lowerval(\Gamma(\max(f,I_O), p')) 
= \lowerval(\Gamma(\max(f,I_{O_{p'}}), p')) \leq g_{i+1}(p') + \epsilon $, 
and for any $w \not\in [H_{i+1}]$, $\max(f,I_O)(w) = f(w) = g_{i+1}(w) $.
Hence by Corollary \ref{Lowervalue-Compare-Corollary},
\begin{equation}
g_i(p) \leq \lowerval(\Gamma(\max(f,I_O), p))
\leq \lowerval(\Gamma_{H_{i+1}}(g_{i+1}, p)) + \epsilon 
\end{equation}
This holds for any $\epsilon>0$, therefore
\begin{equation}
     g_i(p) = \lowerval(\Gamma_{H_{i+1}}(g_{i+1}, p))
\end{equation}
Finally, for any $i \in \N$, and any play $w \not\in [H_i]$, 
we have that $w \not\in [H_{i+1}]$, and hence $g_i(w) = f(w) = g_{i+1}(w)$.
So $\Gamma_{H_i}(g_i)$ is a truncated subgame of $\Gamma_{H_{i+1}}(g_{i+1})$
equivalent for player I.
\\ \ \\
Let $\epsilon>0$.
Since $(\Gamma_{H_i}(g_i))_{i \in \N}$ is a nested series of truncated
subgames equivalent for player I,
by Lemma \ref{Nested-Lowervalue-Lemma} all the games have the same 
lower value, namely $\lowerval(\Gamma_{H_0}(g_0)) = g_0(e) = v$,
and there exists a strategy $\sigma$ for player I
that is $\epsilon$-optimal in all the games $\Gamma_{H_i}(g_i)$, i.e. 
for any strategy $\tau$, and any $i \in \N$,
\begin{equation} 
\expectfin{H_i}{g_i} \geq \lowerval(\Gamma_{H_i}(g_i)) - \epsilon = v-\epsilon
\end{equation}
Now let $\tau$ be any strategy for player II in $\Gamma(\max(f,I_D))$.
Then
\begin{eqnarray}
\lefteqn{ \expectst{\max(f,I_D)} \nonumber} \\
& = & \lim_{i \to \infty} \expectfin{H_i}{h_i} \\
& \geq & \lim_{i \to \infty} \expectfin{H_i}{g_i} \\
& \geq & v - \epsilon
\end{eqnarray}
So $\sigma$ is a strategy for player I of value at least $v-\epsilon$.
This implies that $\lowerval(\Gamma(\max(f,I_D))) \geq v - \epsilon$.
This construction can be done for any $\epsilon>0$, hence
\begin{equation}
\lowerval(\Gamma(\max(f,I_D))) \geq v
\end{equation}
For any $O \supseteq D$, $\lowerval(\Gamma(\max(f,I_D))) \leq
\lowerval(\Gamma(\max(f,I_O)))$, hence
\begin{equation} 
\lowerval(\Gamma(\max(f,I_D))) \leq
\inf_{O \supseteq D, O\mbox{ \scriptsize open}}
\lowerval(\Gamma(\max(f,I_O))) = v 
\end{equation}
Hence $\lowerval(\Gamma(\max(f,I_D))) = v$.
\eop

\begin{cor}
\label{G-Delta-To-Open-Set-Corollary}
Let $S$ be a measurable set, and let $D$ be a $G_\delta$ set.
Suppose that $\Gamma(S \cup D)$ has lower value $v$.
Then for any $\epsilon>0$, there exist an open set $O$, $D \subseteq O$, 
such that $\Gamma(S \cup O)$ has lower value at most $v+\epsilon $.
\end{cor}

\proof
Take $f \equiv I_S$ and apply the non-trivial part of
Theorem \ref{G-Delta-And-Function-Theorem}. \eop

\begin{cor}
\label{G-Delta-Corollary}
Let $D$ be a $G_\delta$ set.
Then $\Gamma(D)$ is determined, and
\begin{equation}
\val(\Gamma(D)) = 
\inf_{O \supseteq D, O\mbox{ \scriptsize open}} \val(\Gamma(O))
\end{equation}
\end{cor}
\proof
For any open set $O \supseteq D$, $\Gamma(O)$ is determined and
$\lowerval(\Gamma(D)) \leq \upperval(\Gamma(D)) \leq \val(\Gamma(O))$.
Applying Theorem \ref{G-Delta-And-Function-Theorem} with $f\equiv 0$
yields the Corollary.
\eop

\subsection{$G_{\delta\sigma}$-sets}

In this subsection, we prove the determinacy of $\Gamma(f)$ in the case
that $f$ is the indicator function of a $G_{\delta\sigma}$ set.
Structurally, this proof is similar to the aforementioned
proof by Davis for $G_{\delta\sigma}$ games of perfect information \cite{Davis}.
We apply the results of the previous subsection
to the $G_\delta$ subsets of a $G_{\delta\sigma}$ set.
Corollary \ref{Open-Union-Corollary} takes the place of the compactness
used in Davis' proof.

\begin{thm}
\label{G-Delta-Sigma-Theorem}
Let $S = \bigcup_i D_i$ be a $G_{\delta\sigma}$ set.
Then $\Gamma(S)$ is determined.
\end{thm}

\prc
Each of the $G_{\delta}$ sets composing the
$G_{\delta\sigma}$ set is enlarged to an open set using Corollary 
\ref{G-Delta-To-Open-Set-Corollary}, 
in such a way that at all times the lower value is not increased by more than
$\epsilon$ (compared to the original game),
where $\epsilon$ is arbitrarily small.
The resulting union of open sets is itself open, and hence determined, and
furthermore Corollary \ref{Open-Union-Corollary}
allows us to conclude that the total increase of
the lower value is still not more than $\epsilon$. 
This means that the upper value
of the original game is also only at most $\epsilon$ more than the lower value. 

Note that, unlike the previous proofs, this proof does not produce an
optimal or $\epsilon$-optimal strategy.
\\ \ \\
\proof
Put $v=\lowerval(\Gamma(S))$.
Let $\epsilon>0$.
Using Corollary \ref{G-Delta-To-Open-Set-Corollary},
we can find inductively open sets 
$O_i \supseteq D_i$ such that for all $j \in \N$, 
\begin{equation}
\lowerval(\Gamma(S \cup \bigcup_{i \leq j+1} O_i)) \leq \lowerval(\Gamma(S \cup
\bigcup_{i \leq j} O_i)) + \epsilon/2^j
\end{equation}
Then for all $j \in \N$,
\begin{equation}
\lowerval(\Gamma(S \cup \bigcup_{i \leq j} O_i)) \leq v+\epsilon
\end{equation}
and hence, for all $j \in \N$, 
\begin{equation}
\val(\Gamma(\bigcup_{i \leq j} O_i)) \leq v+\epsilon
\end{equation}
Then by Corollary \ref{Open-Union-Corollary}, 
\begin{equation}
\val(\Gamma(\bigcup_{i \in \N} O_i)) \leq v+\epsilon
\end{equation}
Since $S = \bigcup_{i \in \N} D_i \subseteq \bigcup_{i \in \N} O_i$, 
\begin{equation}
\upperval(\Gamma(S)) \leq \upperval(\Gamma(\bigcup_{i \in \N} O_i))
= \val(\Gamma(\bigcup_{i \in \N} O_i)) \leq v+\epsilon \end{equation}
This is true for any $\epsilon$, hence
$\upperval(\Gamma(S)) = v = \lowerval(\Gamma(S)) $. \eop

\begin{rem}
The proof of Theorem \ref{G-Delta-Sigma-Theorem} shows that
any $G_{\delta\sigma}$ set (and a fortiori any set of lesser complexity)
can be enlarged to an open set
such that the value of the Blackwell game on that set is not increased
by more than an arbitrarily small amount.
A plausible conjecture is, that this holds for any Borel-measurable set.

This conjecture holds in the case of games of Perfect Information.
Such a game, on a Borel-set $S$, is determined and has value $0$ or $1$.
If it has value $0$ then player II has a winning strategy.
The set of plays that cannot occur if player II uses that strategy,
is an open set, and the game on that set has value $0$ as well.

\end{rem}
\Section{The Axiom of Determinacy for Blackwell Games}

For Games of Perfect Information, there exists the Axiom of
Determinacy, which states that any Game of Perfect Information
with finite choice of moves is determined\footnote{
Formally, AD is an axiom about games with countable choice of moves, whose
payoff function is the indicator function of a set $S \subseteq W$.
But in the case of Games of Perfect Information,
determinacy for games with countable choice of moves is equivalent to
determinacy for games with finite choice of moves,
and determinacy for games with 0-1-valued payoff functions is equivalent to
determinacy for games with arbitrary bounded payoff functions.}.
AD has many interesting consequences, such as the existence of an ultrafilter
on $\aleph_1$, the existence of a complete measure on $\R$, 
the non-existence of a sequence of $\aleph_1$ reals, and the
negation of the Axiom of Choice.
We can formulate an analogue of AD with respect to Blackwell Games,
and look at the consequences of that axiom.
But AD is an axiom about games on all subsets of $W$,
not just on the Borel measurable subsets\footnote{
AD with respect only to games on Borel measurable subsets
(and finite sets $X$ and $Y$) is in fact provable from CAC 
\cite{G-Borel-1,G-Borel-2}.}
An analogous axiom for Blackwell Games should therefore not be limited
to games with measurable payoff functions.
Hence we need to extend the concepts of expectation and value for Blackwell
games.

\begin{defn}
Let $\Gamma(f)$ be a Blackwell Game, where $f$ is bounded but
not necessarily Borel measurable.
Let $\sigma$ and $\tau$ be strategies for players I, II.
$\sigma$ and $\tau$ determine a probability measure 
$\mu_{\sigma, \tau}$ on $W$, induced by setting
\begin{equation}
\mu_{\sigma, \tau}[p]
= P\{ w \mid w\mbox{ hits }p\}
=\prod_{i=1}^n \left(
  \sigma(p_{\mid(i-1)})(x_i) \bullet \tau(p_{\mid(i-1)})(y_i) \right)
\end{equation}
for any position $p = (x_1, y_1, \ldots, x_n, y_n) \in P$.

Instead of the expected income of player I, if she plays according to 
$\sigma$ and player II plays according to $\tau$, we now have the
\em lower \em and \em upper expected income \em :
\begin{eqnarray}
E^{-}(\sigma\mbox{ \rm vs }\tau\mbox{ \rm in }\Gamma(f) ) & = & 
\sup_{g \leq f, g \mbox{ \scriptsize measurable}}
        \int g(w)d\mu_{\sigma, \tau}(w) \\
E^{+}(\sigma\mbox{ \rm vs }\tau\mbox{ \rm in }\Gamma(f) ) & = & 
\inf_{g \geq f, g \mbox{ \scriptsize measurable}}
        \int g(w)d\mu_{\sigma, \tau}(w)
\end{eqnarray}
Lower value and upper value are redefined in the obvious way:
\begin{eqnarray}
\lowerval(\Gamma(f)) & = & \sup_{\sigma}\inf_{\tau}
    E^{-}(\sigma\mbox{ \rm vs }\tau\mbox{ \rm in }\Gamma(f) ) \\
\upperval(\Gamma(f)) & = & \inf_{\tau}\sup_{\sigma}
    E^{+}(\sigma\mbox{ \rm vs }\tau\mbox{ \rm in }\Gamma(f) )
\end{eqnarray}
Note that in the case that $f$ is measurable, 
these definitions reduce to the old definitions.
\end{defn}
\begin{defn}
The \em Axiom of Determinacy for Blackwell Games \em (AD-Bl)
is the statement that for every pair of non-empty finite sets $X,Y$,
and every bounded function $f$ on $W=(X \times Y)^\N$,
the Blackwell Game $\Gamma(f)$ is determined, i.e.
\begin{equation}
\lowerval(\Gamma(f)) = \upperval(\Gamma(f))
\end{equation}
\end{defn}

\begin{thm}
\label{AD-Bl-Measurable-Theorem}
Assuming AD-Bl, it follows that 
all sets of reals are Lebesgue measurable.
\end{thm}
\prc
Let $X$ and $Y$ be the set $\{0,1\}$.
Then we can construct a mapping $\phi: W \to [0, 1]$
such that if either $\sigma$ or $\tau$ is the strategy that assigns 
the $\frac{1}{2}$-$\frac{1}{2}$ probability distribution to every position,
then the measure $\mu_{\sigma, \tau}$ induces the Lebesgue measure on $[0, 1]$.
We can then deduce the measurability of a set $S \subseteq [0,1]$ from
the determinacy of the game $\Gamma(\phi^{-1}[S])$.
\\ \ \\
Some of the consequences of Theorem \ref{AD-Bl-Measurable-Theorem} are,
that AD-Bl is not consistent with AC, and that the consistency of
ZF + AD-Bl cannot be proven in ZFC.
These results are all similar to results for AD.
An open problem is that of the relationship between AD and AD-Bl,
whether AD follows from AD-Bl, or vice versa,
or even whether AD-Bl follows from a stronger version of AD such as 
$\mbox{AD}_\R$.
 From a given game of Perfect Information,
we can easily construct a Blackwell game that is `equivalent',
and assuming AD-Bl we can find 
an $\epsilon$-optimal {\em mixed} strategy for that equivalent Blackwell-game.
However, to derive AD from AD-Bl, we need to have a {\em pure} strategy,
and even though we can interpret any mixed strategy
as a probability distribution on pure strategies,
there is no guarantee that any of these pure strategies will do as well
as the mixed strategy against
\em all \em counterstrategies.
\appendix \clearpage
\Section{Notational Conventions}

\begin{itemize}
\item I, II: the two players
\item $X$, $Y$: the two sets out of which I and II choose their moves
\item $Z$: $X \times Y$
\item $w, w'$: plays
\item $W$: the set of all plays
\item $W_n$: the set of all finite plays of length $n$
\item $S$: a subset of $W$
\item $B$: a basic open subset of $W$
\item $O$: an open subset of $W$
\item $D$: a $G_\delta$ subset of $W$
\item $p, p'$: finite plays or positions
\item $e$: the starting position
\item $P$: the set of all positions
\item $H, H'$: sets of positions, usually stopsets or stopsets-to-be
\item $w \supset p$: the play $w$ hits position $p$
\item $p' \supset p$: the position $p'$ follows $p$
\item $p' \supseteq p$: the position $p'$ follows or is equal to $p$
\item $[p]$: the set of plays hitting position $p$
\item $[H]$: the set of plays hitting any position in $H$
\item $p \cat p', p\cat w$:
the plays obtained by concatenating sequences $p$ and $p'$ or $w$
\item $p_{\mid n}, w_{\mid n}$: the plays consisting of the first 
$n$ moves in $p$ or $w$
\item $\len(p)$: the length of a position $p$
\\
\item $f, g, h$: payoff functions
\item $f(w)$: the payoff player I gets from player II for the play $w$
\item $f(p)$: if $p$ is a stopping position, the payoff player I gets from 
player II for any play $w$ that passes through $p$
\item $\Gamma(f)$: the game with payoff function $f$
\item $\Gamma_H(f)$: the game with payoff function $f$ and stopset $H$
\item $\Gamma_n(f)$: the finite game with payoff function $f$ and 
stopset $W_n$
\item $\Gamma(S)$: the game with payoff function $I_S$
\item $\Gamma(f, p)$: the game with payoff function $f$, starting from
position $p$
\item $\sigma$: a strategy for player I
\item $\tau$: a strategy for player II
\item $\expectst{f}$: the expected payoff when playing $\sigma$ 
against $\tau$ in the game $\Gamma(f)$
\item $\val(\sigma \mbox{ in }\Gamma(f))$: the value of $\sigma$ 
in the game $\Gamma(f)$
\item $\val(\Gamma(f))$: the value of the game $\Gamma(f)$
\item $\lowerval(\Gamma(f))$: the lower value of the game $\Gamma(f)$
\item $\upperval(\Gamma(f))$: the upper value of the game $\Gamma(f)$
\item $u, v$: values
\end{itemize}

\Section{Proofs}
\subsection{Equivalent Truncated Subgames}
\begin{trivlist} \item
{\bfseries Lemma \ref{Lowervalue-Lemma}}
{\itshape
Let $\Gamma(f)$ be a Blackwell game,
and let $\Gamma_H(g)$ be a truncated subgame of $\Gamma(f)$,
truncated at a set of positions $H$, equivalent for player I [for player II].
Then $\lowerval(\Gamma(f))=\lowerval(\Gamma_H(g))$ 
[$\upperval(\Gamma(f))=\upperval(\Gamma_H(g))$]. Furthermore,
for any $\epsilon>0$, 
any $\epsilon$-optimal strategy for player I [player II] in $\Gamma_H(g)$
(if it is undefined on all positions at or after positions in $H$) 
can be extended to an $\epsilon$-optimal strategy for player I [player II] 
in $\Gamma(f)$.
}
\end{trivlist}
\proof
Let $\Gamma(f)$ be a Blackwell game,
and let $\Gamma_H(g)$ be a truncated subgame of $\Gamma(f)$,
truncated at a set of positions $H$, equivalent for player I.
This means that for any $p \in H$, $g(p)=\lowerval(\Gamma(f,p))$, and
for any $w \not\in [H]$, $g(w)=f(w)$. 

Let $\epsilon>0$, and let $\sigma_0$ be an $\epsilon$-optimal strategy
for player I in the game $\Gamma_H(g)$. If $v=\lowerval(\Gamma_H(g))$, and
$u = \val(\sigma_0\mbox{ in }\Gamma_H(g))$, then $0 \leq v-u < \epsilon$.
So choose $\delta>0$ such that $\delta < \epsilon-(v-u)$, 
and fix for each $p \in H$ a
$\delta$-optimal strategy $\sigma_p$ for player I in $\Gamma(f,p)$.
By Remark \ref{Undefined-Strategy-Remark}
we may assume that for any $p \in H$, 
$\sigma_p$ is defined exactly on $p$ and those positions that are after $p$, 
and that $\sigma_0$ is defined exactly on those positions that are
not at or after \em any \em position in $H$.
It follows that $\sigma=\bigcup_{p\in H}\sigma_p \cup \sigma_0$ is
a well-defined strategy for player I in the game $\Gamma(f)$.

Now let $\tau$ be a strategy for player II in $\Gamma(f)$.
Then for each $p \in H$,
\begin{eqnarray}
\int_{w\in[p]} f(w) d\mu_{\sigma, \tau}(w)
& = & \left( \frac{\int_{w\in[p]} f(w) d\mu_{\sigma, \tau}(w)}
                    {\mu_{\sigma, \tau}[p]} \right)
                     \mu_{\sigma, \tau}[p] \\
& = & \left( \int_{w\in[p]} f(w)
              d\mu_{\sigma, \tau\mbox{ \scriptsize in }\Gamma(f,p)}(w) \right)
                    \mu_{\sigma, \tau}[p] \heq \\
& = & \expectst{f, p}
                    \mu_{\sigma, \tau}[p] \\
& = & \expect{\sigma_p}{\tau}{\Gamma(f, p)}
                    \mu_{\sigma_0, \tau}[p] \\
& \geq & \left( \lowerval(\Gamma(f, p)) - \delta \right)
                    \mu_{\sigma_0, \tau}[p] \\
& = & (g(p)-\delta) \mu_{\sigma_0, \tau}[p] 
\end{eqnarray}
and consequently,
\begin{eqnarray}
\lefteqn{\expectst{f}} \nonumber \\
& = & \int_{w\in W} f(w)d\mu_{\sigma, \tau}(w) \\
& = & \sum_{p \in H}
        \left( \int_{w\in[p]} f(w) d\mu_{\sigma, \tau}(w) \right)
      + \int_{w\not\in[H]} f(w) d\mu_{\sigma, \tau}(w) \\
& \geq & \sum_{p \in H} (g(p)-\delta)\mu_{\sigma_0, \tau}[p]
      + \int_{w\not\in[H]} g(w) d\mu_{\sigma_0, \tau}(w) \\
& \geq & \left( \sum_{p \in H}g(p)\mu_{\sigma_0, \tau}[p]
  + \int_{w\not\in[H]} g(w) d\mu_{\sigma_0, \tau}(w) \right) - \delta \heq \\
& = & \int_w g(w)d\mu_{\sigma_0, \tau}(w) - \delta \\
& \geq & u-\delta \\
& > & v-\epsilon 
\end{eqnarray}
So $\sigma$ is an extension of $\sigma_0$ of value greater than $v-\epsilon$.
Since this can be done for any $\epsilon>0$, and any $\epsilon$-optimal
strategy $\sigma_0$, this implies that 
$\lowerval(\Gamma(f)) \geq v$.

Similarly, let $\epsilon>0$, and let $\sigma$ be a strategy for player I
in $\Gamma(f)$.
Then we can find counterstrategies $\tau_0$, $\tau_p$ for $p \in H$ for
player II in $\Gamma_H(g)$, $\Gamma(f, p)$ respectively, such that
$\expect{\sigma}{\tau_0}{\Gamma_H(g)} < \lowerval(\Gamma_H(g))+\epsilon/2$, and
$\expect{\sigma}{\tau_p}{\Gamma(f, p)} < \lowerval(\Gamma(f, p))+\epsilon/2$
for any $p \in H$.
We can assume that for any $p \in H$, 
$\tau_p$ is defined exactly on those positions that are at or after $p$, 
and that $\tau_0$ is defined exactly on those positions that are
not at or after any \em any \em position in $H$.
It follows that $\tau=\bigcup_{p\in H}\tau_p \cup \tau_0$ is
a well-defined strategy for player II in the game $\Gamma(f)$.
Then for each $p \in H$,
\begin{eqnarray}
\int_{w\in[p]} f(w) d\mu_{\sigma, \tau}(w)
& = & \left( \frac{\int_{w\in[p]} f(w) d\mu_{\sigma, \tau}(w)}
                    {\mu_{\sigma, \tau}[p]} \right)
                     \mu_{\sigma, \tau}[p] \\
& = & \left( \int_{w\in[p]} f(w)
                 d\mu_{\sigma, \tau\mbox{ \scriptsize in }\Gamma(f,p)}(w) \right)
                    \mu_{\sigma, \tau}[p] \heq \\
& = & \expectst{f, p}
                    \mu_{\sigma, \tau}[p] \\
& = & \expect{\sigma}{\tau_p}{\Gamma(f, p)}
                    \mu_{\sigma, \tau_0}[p] \\
& \leq & \left( \lowerval(\Gamma(f, p)) + \epsilon/2 \right)
                    \mu_{\sigma, \tau_0}[p] \\
& = & (g(p)+\epsilon/2)\mu_{\sigma, \tau_0}[p]
\end{eqnarray}
and consequently,
\begin{eqnarray}
\lefteqn{\expectst{f}} \nonumber \\
& = & \int_{w\in W} f(w)d\mu_{\sigma, \tau}(w) \heq \\
& = & \sum_{p \in H}
        \left( \int_{w\in[p]} f(w) d\mu_{\sigma, \tau}(w) \right)
      + \int_{w\not\in[H]} f(w) d\mu_{\sigma, \tau}(w) \\
& \leq & \sum_{p \in H} (g(p)+\epsilon/2)\mu_{\sigma, \tau_0}[p]
      + \int_{w\not\in[H]} g(w) d\mu_{\sigma, \tau_0}(w) \\
& \leq & \sum_{p \in H}g(p)\mu_{\sigma, \tau_0}[p]
      + \int_{w\not\in[H]} g(w) d\mu_{\sigma, \tau_0}(w) + \epsilon/2 \heq \\
& = & \int_w g(p)d\mu_{\sigma, \tau_0}(w) +\epsilon/2 \\
& < & v+\epsilon/2+\epsilon/2 \\
& = & v+\epsilon 
\end{eqnarray}
Since this can be done for any strategy $\sigma$ and any $\epsilon>0$,
it follows that $\lowerval(\Gamma(f)) \leq v$.
Hence
\begin{equation}
\lowerval(\Gamma(f)) = v=\lowerval(\Gamma_H(g))
\end{equation}
An $\epsilon$-optimal strategy for player I in $\Gamma_H(g)$ can be extended
to an $\epsilon$-optimal strategy for player I in $\Gamma(f)$ as in the
first part of the proof.

The proof for the upper value is analogous.
\eop
\begin{trivlist} \item
{\bfseries Theorem \ref{Nested-Uppervalue-Lemma}}
{\itshape
Let $(\Gamma_{H_i}(g_i))_{i \in \N}$ be a nested series of truncated games
equivalent for player I [player II]. 
Then all the games have the same lower value [upper value].
Furthermore, we can find a strategy for player I [player II] 
that is $\epsilon$-optimal in all the games $\Gamma_{H_i}(g_i)$.
}
\end{trivlist}
\proof
For all $i \in \N$, $\Gamma_{H_i}(g_i)$ is a truncated subgame of
$\Gamma_{H_{i+1}}(g_{i+1})$, equivalent for player I.
By Lemma \ref{Lowervalue-Lemma}, this means that for all $i \in \N$,
$\lowerval(\Gamma_{H_i}(g_i))=\lowerval(\Gamma_{H_{i+1}}(g_{i+1}))$.
It follows that for all $i \in \N$,
\begin{equation}
\lowerval(\Gamma_{H_i}(g_i))=\lowerval(\Gamma_{H_0}(g_0))
\end{equation}
Now let $\epsilon>0$, and let $\sigma_0$ be an $\epsilon$-optimal strategy
for player I in $\Gamma_{H_0}(g_0)$.
Assume that $\sigma_0$ is only defined on non-stopping positions of 
$\Gamma_{H_0}(g_0)$. 
We can inductively define $\epsilon$-optimal
strategies $\sigma_i$ in the games $\Gamma_{H_i}(g_i)$ such that
for all $i \in \N$, $\sigma_{i+1}$ is an extension of $\sigma_i$.
For suppose we have extended $\sigma_0$ to an $\epsilon$-optimal strategy 
$\sigma_i$ in $\Gamma_{H_i}(g_i)$, defined only on non-stopping positions.
Then $\sigma_i$ is undefined on positions at or after positions in $H_i$,
and thus by Lemma \ref{Lowervalue-Lemma} we can extend $\sigma_i$ to an 
$\epsilon$-optimal strategy $\sigma'_{i+1}$ in $\Gamma_{H_{i+1}}(g_{i+1})$.

By Remark \ref{Undefined-Strategy-Remark}, we may restrict this strategy
$\sigma'_{i+1}$ to the non-stopping positions in $\Gamma_{H_{i+1}}(g_{i+1})$.
This restriction $\sigma_{i+1}$ is an extension of $\sigma_i$ as well.
For let $p$ be a stopping position in $\Gamma_{H_{i+1}}(g_{i+1})$.
If $p$ is preceded by a stopping position in $\Gamma_{H_i}(g_i)$, then
$p$ is itself a stopping position in $\Gamma_{H_i}(g_i)$.
Otherwise, we observe that for any play $w \supset p$, $g_{i+1}(w)=g_{i+1}(p)$,
and for any position $q \supseteq p$,
$\Gamma(g_{i+1},q)$ is determined with value $g_{i+1}(p)$.
By the truncated-subgame condition, this implies that 
any $w \in [p]$ has payoff $g_{i+1}(p)$ in $\Gamma_{H_i}(g_i)$, 
and again $p$ is a stopping position.

Now, if $p$ is a position, and $\sigma_i(p)$, $\sigma_j(p)$ are defined,
$i<j$, then $\sigma_i(p) = \sigma_{i+1}(p) = \ldots = \sigma_j(p) $.
Therefore we can define $\sigma$ by $\sigma=\bigcup_{i \in \N} \sigma_i$, i.e.
\begin{equation}
\sigma(p) = \sigma_i(p)
     \mbox{ if $\sigma_i$ is defined on $p$, for any $i \in \N$ }
\end{equation}
Then for any $i \in \N$, for any non-stopping position $p$ of $\Gamma_{H_i}(g_i)$,
$\sigma(p)$ is defined and $\sigma(p) = \sigma_i(p)$. Hence, $\sigma$ is
$\epsilon$-optimal in all the games. $\sigma$ is defined on all the 
non-stopping positions of all the games $\Gamma_{H_i}(g_i)$, but if necessary 
we can extend $\sigma$ to a strategy $\sigma'$ defined on all positions,
by choosing a default probability distribution for those positions on which 
$\sigma$ is not defined.
\eop
\subsection{Finite Games}
\begin{lem}
\label{Convex-Separation}
Let $C$ be a closed convex set in $\R^n$, and let $b \in \R^n-C$.
Then there exists a hyperplane separating $b$ and $C$, i.e. a vector 
$y \in \R^n-\{0\}$ and $d \in \R$ such that 
\begin{eqnarray}
\label{convexeq1}
y^Tb & > & d \\
\label{convexeq2}
y^Tz & < & d \mbox{ for each }z \in C
\end{eqnarray}
\end{lem}
\proof
Since the theorem is trivial if $C = \emptyset$, we assume $C \neq \emptyset$.
As $C$ is closed, there exists a vector $c$ in $C$ 
that is nearest to $b$, i.e. that minimizes $\|z-b\|$. 
And because $b \not\in C$, $\|c-b\| > 0$.
We define
\begin{eqnarray}
y & = & b-c \\
d & = & \frac{1}{2}(y^Tb+y^Tc) 
\end{eqnarray}
Then $y^Tb \;-d = y^T(b-\frac{1}{2}(b+c)) = \frac{1}{2}(b-c)^T(b-c) > 0$, 
proving that (\ref{convexeq1}) holds.

Now suppose that (\ref{convexeq2}) does not hold.
Then for some $z \in C$, $y^Tz \geq d$. Since,
$y^Tc \;-d = y^T(c-\frac{1}{2}(b+c)) = -\frac{1}{2}(b-c)^T(b-c) < 0$,
this implies $y^T(z-c)>0$.
Hence there exists a $\lambda$ with $0 < \lambda \leq 1$ and
\begin{equation}
\lambda < \frac{2y^T(z-c)}{\|z-c\|^2}
\end{equation}
As C is convex, $c+\lambda(z-c)$ belongs to $C$. Moreover,
\begin{eqnarray}
\|(c+\lambda(z-c)) - b\|^2 & = & \|\lambda(z-c)-y\|^2 \\ 
& = & \lambda^2\|z-c\|^2-2\lambda y^T(z-c)+\|y\|^2  \\
& < & \|y\|^2 = \|c-b\|^2
\end{eqnarray}
contradicting the fact that $c$ is a point in $C$ nearest to $b$.
Therefore, for all $z \in C$, $y^Tz < d$, and (\ref{convexeq2}) holds.
\eop
\begin{trivlist} \item
{\bfseries Theorem \ref{One-Move} (Von Neumann's Minimax Theorem\cite{Neumann})}
{\itshape
Let $\Gamma_1(f)$ be a finite one-round Blackwell game (i.e. of length 1).
Then $\Gamma_1(f)$ is determined, and both players have optimal strategies.
}
\end{trivlist}
\proof
$f$ is (or can be interpreted as) a function $X \times Y \to \R$.
Without loss of generality we may assume that $X = \{1, \ldots, n\},
Y = \{1, \ldots, m\}$. 
A strategy $\sigma$ for player I in $\Gamma_1(f)$ is (or can be interpreted as)
a nonnegative vector $(x_1, \ldots, x_n)$ such that $\sum_{i=1}^n x_i =1$,
and similarly, a strategy $\tau$ for player II is (or can be interpreted as)
a nonnegative vector $(y_1, \ldots, y_m)$ such that $\sum_{j=1}^m y_j =1$.
It is easily seen that
\begin{eqnarray}
\expectfin{1}{f} & = & 
   \sum_{i=1}^n\sum_{j=1}^m x_i y_j f(i,j) \\
\val(\sigma\mbox{ in }\Gamma_1(f)) & = &
   \min_{1 \leq j \leq m} \sum_{i=1}^n x_i f(i, j) \\
\val(\tau\mbox{ in }\Gamma_1(f)) & = &
   \max_{1 \leq i \leq n} \sum_{j=1}^m y_j f(i, j)
\end{eqnarray}
Let $C \subset \R^m$ be the convex hull of the vectors 
$(f(i,1), \ldots, f(i,m))$, $i = 1, \ldots, n$. Then each point $z \in C$
corresponds to (at least) one strategy $\sigma$ for player I, 
with $\val(\sigma\mbox{ in }\Gamma_1(f))$ being equal to the
coordinate of $z$ which has the least value.
Let $C^{-} = \{ z \in \R^m \mid \exists z' \in C : z \leq z' \}$,
using coordinatewise ordening of vectors.
Then player I has a strategy of value $u \in \R$ iff 
the vector $(u, \ldots, u) \in C^{-}$.
It is obvious that $C$ and $C^{-}$ are both closed convex sets.

Now let $v = \lowerval(\Gamma_1(f))$, $\epsilon>0$.
Then there exists no strategy $\sigma$ for player I of value $v+\epsilon$, 
i.e. $b = (v+\epsilon, \ldots, v+\epsilon) \not\in C^{-}$.
By Lemma \ref{Convex-Separation} this implies that 
there exist a vector $y \in \R^m-\{0\}$ and $d \in \R$, such that 
$y^Tb > d$ while $y^Tz < d$ for any $z \in C^{-}$.
For any $j \leq m$, any $z \in C^{-}$, and any $r>0$,
$y^Tz - y_jr = y^T(z_1, \ldots, z_{j-1}, z_j-r, z_{j+1}, \ldots, z_m) < d$,
i.e. $y_j > (y^Tz - d)/r$.
Since $r>0$ is arbitrary, it follows that for any $j \leq m$, $y_j \geq 0$.
Since $y \neq 0$, we may also assume, without loss of generality, 
that $\sum_{j=1}^m y_j = 1$.
Then for any $z \in C^{-}$, $y^Tz < d < y^Tb = v+\epsilon$.
In particular, for $i=1, \ldots, n$, $\sum_{j=1}^m y_j f(i,j) < v+\epsilon$.
Thus, $y$ corresponds to a strategy $\tau$ of value lower than $v+\epsilon$.
Since we can do this construction for any $\epsilon>0$, it
follows that $\upperval(\Gamma_1(f)) = v = \lowerval (\Gamma_1(f))$.

The existence of an optimal strategy for player I follows from the observation
that $C$ is a closed and bounded subset of $\R^m$ (and hence compact), and that
the function $\min: \R^m \to \R$ taking the minimum coordinate is continuous.
A similar argument yields the existence of an optimal strategy for player II.
\eop
\begin{trivlist} \item
{\bfseries Theorem \ref{Finite-Theorem}}
{\itshape
Let $\Gamma_n(f)$ be a finite Blackwell game of length $n$.
Then $\Gamma_n(f)$ is determined, and both players have optimal strategies.
}
\end{trivlist}
\proof
We will proof this by induction on the length $n$ of the games.

For $n=0$, determinacy is trivial.
So let $n>0$, let $\Gamma_n(f)$ be a finite Blackwell game of length $n$,
and suppose that all finite Blackwell games of length $m < n$
have already been shown to be determined, and to have optimal strategies.
For each position $p \in W_1$, the game $\Gamma(f, p)$ is finite of 
length $\leq n-1$. Hence, by the induction hypothesis, each of those games is
determined and has a value $\val(\Gamma(f, p))$. Define the payoff function
$g: W \to \R$ by
\begin{equation}
g(p) = \val(\Gamma(f, p))\mbox{ for }p \in W_1
\end{equation}
Then by Von Neumann's Minimax Theorem, the game $\Gamma_1(g)$ is determined.
$\Gamma_1(g)$ is by its definition an equivalent truncated subgame of
$\Gamma_n(f)$,
so by Corollary \ref{Value-Corollary}, $\Gamma_n(f)$ is determined.
Furthermore, by Remark \ref{Strategy-Extension-Optimal-Remark}
the optimal strategy produced by Von Neumann's Mini-Max Theorem can be extended
to an optimal strategy in $\Gamma_n(f)$.
\eop

\subsection{Lebesgue Measurability from AD-Bl}
\begin{trivlist} \item
{\bfseries Theorem \ref{AD-Bl-Measurable-Theorem}}
{\itshape
Assuming AD-Bl, it follows that 
all sets of reals are Lebesgue measurable.
}
\end{trivlist}
\proof
It suffices to show that the Lebesgue measure on $[0, 1]$ is complete.
Set $X = Y = \{0,1\}$, and define $\phi: W \to [0, 1]$ by
\begin{equation}
\phi((x_1, y_1, x_2, y_2, \ldots)) = \sum_{i=1}^\infty 2^{-i}(x_i \oplus y_i)
\end{equation}
where $0 \oplus 0 = 1 \oplus 1 = 0$, $0 \oplus 1 = 1 \oplus 0 = 1$.
Now let $\sigma, \tau$ be strategies, and suppose that one of those strategies
is the strategy that assigns the $\frac{1}{2}$-$\frac{1}{2}$ probability
distribution on $X$ or $Y$, respectively.
Then for any $i \in \N$, $x_i \oplus y_i$ has equal chances of being $0$ or $1$,
and the distribution of $\phi(w)$ on $[0, 1]$ is the uniform
distribution on $[0, 1]$ under the Lebesgue measure on $[0, 1]$.
It follows that for any subset $S \subset [0, 1]$,
\begin{eqnarray}
\mu^{\mbox{\scriptsize inner}}{\sigma, \tau}(\phi^{-1}[S])
& = & \mu^{\mbox{\scriptsize inner}}_{\lambda}(S) \\
\mu^{\mbox{\scriptsize outer}}{\sigma, \tau}(\phi^{-1}[S])
& = & \mu^{\mbox{\scriptsize outer}}_{\lambda}(S)
\end{eqnarray}
where $\mu^{\mbox{\scriptsize inner}}(A)=
\sup_{B \subseteq A \mbox{ \scriptsize measurable}}\mu(B)$,
$\mu^{\mbox{\scriptsize outer}}(A)=
\inf_{B \supseteq A \mbox{ \scriptsize measurable}}\mu(B)$.

Let $S \subset [0, 1]$.
No strategy for player I in the game $\Gamma(\phi^{-1}[S])$
can have value greater than $\mu^{\mbox{\scriptsize inner}}_{\lambda}(S)$,
since this is the lower expected income of any strategy for player I
against the $\frac{1}{2}$-$\frac{1}{2}$ strategy.
Similarly, no strategy for player II in the game $\Gamma(\phi^{-1}[S])$
can have value less than $\mu^{\mbox{\scriptsize outer}}_{\lambda}(S)$.
Therefore
\begin{equation}
\lowerval(\Gamma(\phi^{-1}[S]))
\leq \mu_{\lambda}^{\mbox{\scriptsize inner}}(S)
\leq \mu_{\lambda}^{\mbox{\scriptsize outer}}(S)
\leq \upperval(\Gamma(\phi^{-1}[S])).
\end{equation}
 From the determinacy of the game $\Gamma(\phi^{-1}[S])$,
it now follows that
\begin{equation}
\mu_{\lambda}^{\mbox{\scriptsize inner}}(S) = \mu_{\lambda}^{\mbox{\scriptsize outer}}(S)
\end{equation}
Since this holds for arbitrary sets $S \subset [0, 1]$,
all subsets of $[0, 1]$ are Lebesgue measurable.
\eop 
\Section{The Locally Optimal Strategy}

As can be seen by inspection of the proof of Theorem \ref{Finite-Theorem},
an optimal strategy in a finite game is to use, at each position $p$,
Von Neumann's Minimax Theorem to generate the probability distribution for
that position $p$. 
We do this by calculating, for each position $p$, the optimal strategy
for the one-round game whose payoffs are the values of the games
starting at the positions directly following $p$.

In a non-finite game, we can do the same thing,
i.e. for a given position $p$, we can apply Von Neumann's Minimax Theorem
to the one-round game whose payoffs are the lowervalues of the games starting
 from the positions directly following $p$.
By Lemma \ref{Lowervalue-Lemma}, the value of this one-round game is equal to
the lowervalue of the game starting at $p$, giving us a relationship
between lowervalues at different positions.
In an open game, for player II, the strategy of \em always \em using
this `locally optimal' probability distribution is an optimal strategy
(Remark \ref{Open-Local-Strategy-Remark}).
But this does not hold in general: 
see also Example \ref{Anti-Local-Strategy-Example} for a counterexample.

As to why this `local optimality' does not translate to optimality for
the whole strategy:
a single probability distribution
that is optimal in the corresponding one-round game
would be optimal in the whole game,
if for any position $p$ following the current position,
the strategy used in the subgame $\Gamma(f, p)$ were optimal.
But in order to use this to show that each of the separate probability
distributions in the `locally optimal strategy' is optimal in the whole game,
we would first need that
the strategy as a whole is optimal in entire subgames,
and for that we need that
each of the separate probability distributions is optimal in those subgames,
and for that we need that
the strategy as a whole is optimal in entire subgames of subgames,
etcetera, etcetera, ad infinitum.
In the case of finite games this eventually reduces to subgames of length 0,
which are trivial, but in the case of infinite games the proof falls through.

In general, it is impossile to construct an optimal or $\epsilon$-optimal
strategy merely from the values a game has at each of its positions,
as the following example shows.
\begin{ex}
\label{Complementary-Example}
Consider the following two Blackwell games.
Player II generates a sequence of 0's and 1's.
In the game $\Gamma(S_1)$, player II wins if he generates infinitely many 1's.
In the game $\Gamma(S_2)$, player II wins if he generates only finitely many 1's.
Player I has no influence over the outcome of the game.
Payoff is $1$ if player I wins, $0$ if player II wins.

It is clear that player II can win from any position, in both games.
Hence both games have lower and upper value 0, starting from all positions.
But any strategy for player II that is good in $\Gamma(S_1)$, will be bad in
$\Gamma(S_2)$, and vice versa, since the two sets of winning positions
are complementary.
Hence there cannot be any method of finding optimal or good strategies that
merely uses the values of a game.

Note that here we have two sets of winning positions $S_1$ and $S_2$,
such that $\Gamma(S_1)$, $\Gamma(S_2)$ have value $0$ (in all positions),
and $\Gamma(S_1 \cup S_2)$ has value $1$ (in all positions).
\end{ex}


\begin{thebibliography}{10}
\pagestyle{myheadings} 
\thispagestyle{myheadings} 
\markboth{\sc Bibliography}{\sc Blackwell Games}

\bibitem{Connect-Four}
{\bf V. Allis},
{\em A Knowledge-based Approach of Connect-Four,}
Report IR-163,
Faculty of Mathematics and Computer Science
at the Vrije Universiteit Amsterdam, The Netherlands,
1988.

\bibitem{Blackwell}
{\bf D. Blackwell}, {\em Infinite $G_\delta$ games with imperfect information},
Zastosowania Matematyki Applicationes Mathematicae,
Hugo Steinhaus Jubilee Volume X, p. 99-101,
1969.

\bibitem{Blackwell-2}
{\bf D. Blackwell},
{\em Operator Solution of Infinite $G_\delta$ Games of Imperfect Information},
in: T.W. Anderson, K.B. Athreya, D.L. Iglehart (editors),
{\em Probability, Statistics and Mathematics. Papers in Honor of S. Karlin},
p. 83-87, Academic Press, New York,
1989.

\bibitem{Blackwell-3}
{\bf D. Blackwell} and {\bf M.A. Girshick},
{\em Theory of Games and Statistical Decisions},
Dover Publications Inc., New York,
1954.

\bibitem{Davis}
{\bf M. Davis},
{\em Infinite games of perfect information},
in: M. Dresher, L.S. Shapley, A.W. Tucker (editors),
{\em Advances in game theory},
Annals of Mathematics Studies 52, p.85-101,
Princeton University Press, Princeton, N.J.,
1964.

\bibitem{G-Borel-1}
{\bf A.J.C. Hurkens}, {\em Borel determinacy without the axiom of choice},
Universiteitsdrukkerij Nijmegen, Nijmegen, The Netherlands,
1993.

\bibitem{G-Borel-2}
{\bf D.A. Martin}, {\em Borel Determinacy},
Annals of Mathematics 102, p. 363-371,
1975.

\bibitem{G-Borel-3}
{\bf D.A. Martin}, {\em A purely inductive proof of Borel determinacy},
in: A. Nerode and R.A. Shore (Editors), {\em Recursion Theory},
Proceedings of Symposia in Pure Mathematics 42, p. 303-308,
American Mathematical Society, Providence, R.I.,
1985.

\bibitem{AD-1}
{\bf D.A. Martin} and {\bf J.R. Steel},
{\em A proof of projective determinacy},
Journal of the American Mathematical Society 2, p. 71-125,
1989.

\bibitem{AxiomOfChoice}
{\bf D. Gale} and {\bf F.M. Stewart}, {\em Infinite Games with Perfect
Information}, in: {\em Contributions to the theory of Games},
Annals of Mathematics Studies 28, p.245-266,
Princeton University Press, Princeton, N.J.,
1953.

\bibitem{AD-2}
{\bf J. Mycielski} and {\bf H. Steinhaus},
{\em A mathematical axiom contradicting the axiom of choice},
Bulletin de l'Acad\'{e}mie Polonaise des Sciences, 
S\'{e}rie des sciences math\'{e}matiques, astronomiques et physiques 10,
p. 1-3,
1962.

\bibitem{Neumann}
{\bf J. von Neumann}, {\em Zur Theorie der Gesellschaftsspiele},
Math. Annalen, Volume 100, p. 295-320,
1928.

\end{thebibliography}
\end{document}